\documentclass[12pt,thmsa,emstex]{article}
\UseRawInputEncoding
\usepackage{amsfonts}
\usepackage{t1enc}
\usepackage{amssymb}
\usepackage{epsfig}
\usepackage[french,english]{babel}
\usepackage{amsmath}
\usepackage{graphics}
\usepackage{layout}
\usepackage{latexsym}
\usepackage{color}
\usepackage{verbatim}
\usepackage{authblk}

\def\mylabel#1{\label{#1}}

\newtheorem{theorem}{Theorem} 
\newtheorem{lemma}[theorem]{Lemma}
\newtheorem{corollary}[theorem]{Corollary}
\newtheorem{proposition}[theorem]{Proposition}

\newtheorem{exercise}[theorem]{Exercise}

\newtheorem{remark}{Remark}
\newtheorem{example}{\bf{Example}}

\def\bit{\begin{itemize}}
\def\eit{\end{itemize}}

\reversemarginpar   

\def\bc{\begin{center}}
\def\ec{\end{center}}

\def\bthm{\begin{theorem}}
\def\ethm{\end{theorem}}

\def\bcor{\begin{corollary}}
\def\ecor{\end{corollary}}

\def\bprop{\begin{proposition}}
\def\eprop{\end{proposition}}

\def\blem{\begin{lemma}}
\def\elem{\end{lemma}}

\def\bex{\begin{example}}
\def\eex{\end{example}}

\def\bexo{\begin{exercise}}
\def\eexo{\end{exercise} }

\def\brem{\begin{remark}}
\def\erem{\end{remark}}
\def\prf{{\bf Proof: }}

\def\bdes{\begin{description}}
\def\edes{\end{description}}

\def\ita{\item[(a)]}
\def\itb{\item[(b)]}

\def\iti{\item[(i)]}
\def\itii{\item[(ii)]}

\def\beq{\begin{equation}}
\def\eeq{\end{equation}}

\def\ben{\begin{enumerate}}
\def\een{\end{enumerate}}

\def\beqar{\begin{eqnarray}}
\def\eeqar{\end{eqnarray}}

\def\beqarr{\begin{eqnarray*}}
\def\eeqarr{\end{eqnarray*}}

\def\qed{\hfill $\Box$ \\[2ex]}
\def\prf{{\bf Proof: }\hspace{.1in}}

\catcode`\@=11

\def\ZZ{{\mathbb Z}}       
\def\RR{{\mathbb R}}  
\def\Rp{{\mathbb R}_+}   

\def\NN{{\mathbb N}}

\def\rar{\rightarrow}
\def\eps{\varepsilon}

\begin{document}
\title{A note on the top Lyapunov exponent of linear cooperative systems}
\author[1]{\small Michel Bena{\"i}m }
\author[2]{Claude Lobry}
\author[3]{Tewfik Sari}
\author[4]{\'Edouard Strickler}
\affil[1]{Institut de Math\'ematiques, Universit\'e de Neuch{\^a}tel, Switzerland}
\affil[2]{C.R.H.I, Universit\'e Nice Sophia Antipolis, France}
\affil[3]{ITAP, University of Montpellier, INRAE, Institut Agro, Montpellier, France}
\affil[4]{ Universit\'e de Lorraine, CNRS, Inria, IECL, Nancy, France }

{\small \date{\today}}

\maketitle
\abstract{In a recent paper \cite{carmona2022},  P. Carmona gives an asymptotic formulae for the top Lyapunov exponent of a linear $T$-periodic cooperative differential equation, in the limit $T \rar \infty.$ This short note discusses and extends this result.
The assumption that the system is $T$-periodic is replaced by the more general assumption that it is driven by a continuous time uniquely ergodic Feller Markov process $(\omega_{t})_{t > 0}.$ When $\omega_{t}$ is replaced by $\omega^T_{t} = \omega_{t/T},$ asymptotic formulas  for the top Lyapunov exponent in the fast (i.e~$T \rar \infty$) and slow ($T \rar 0$) regimes are given.} 
\section{Notation and main results}
 Let $d \geq 1$ be an integer. Let ${\cal M}$ denote the closed convex cone consisting of real $d \times d$ matrices having off diagonal nonnegative entries. Elements of ${\cal M}$ are usually called {\em Metzler} matrices. As usual, a matrix $M \in {\cal M}$ is called {\em irreducible} if for all $i,j \in \{1,\ldots,d\}$ there exist $n \in \NN$ and
 a sequence $i_1 = i, i_2, \ldots, i_n = j$ such that $M_{i_l, i_{l+1}} > 0$ for $l = 1, \ldots, n-1.$ Equivalently $e^M$ has positive entries.
 Throughout, we let $S$ denote a compact metric space and
$$A: S \rar {\cal M},$$
  a continuous mapping.
  We consider the linear differential equation
  \beq
  \label{eq:ode}
     \frac{dy}{dt}=A(\omega_{t}) y
  \eeq
  with initial condition $y(0) = x  \in \RR^d_+ \setminus \{0\},$ under the following  assumptions:
  \bdes
  \iti The process $(\omega_t)_{t \geq 0}$ is a  continuous time  Feller Markov process\footnote{The precise definition will be recalled in the beginning of Section \ref{sec:proofs}}  on $S$  and is
  {\em uniquely ergodic}. By this, we mean that $(\omega_t)_{t \geq 0}$ has a unique invariant probability measure denoted $\mu.$
  \itii The average matrix $\bar{A} = \int_S A(s) \mu(ds)$ is irreducible.
  \edes
  \brem
  {\rm A sufficient (but non necessary) condition ensuring that $\bar{A}$ is irreducible is that $A(s)$ is irreducible for some $s$ in the topological support of $\mu.$ The (easy) proof is left to the reader.}
  \erem
The assumption that $A(s)$ is Metzler for all $s \in S,$ makes the non-autonomous differential equation (\ref{eq:ode}) {\em cooperative} in the sense that $\frac{\partial \dot{y}_i}{\partial y_j} \geq 0$ for all $i \neq j$ (we refer the reader to \cite{HirschSmith06} for a comprehensive introduction to the theory of deterministic cooperative systems). Systems of this form naturally occur in population dynamics where individuals can migrate between different patches (see e.g~\cite{BLSS23a}, \cite{BLSS23b} and references therein) or different states (see e.g~ \cite{MalikSmith08}, \cite{carmona2022}). They also occur as linearized systems of non-linear cooperative systems (for instance in certain  epidemic models \cite{BS2019}). In all these settings, the process $(\omega_t)_{t \geq 0}$  represents the time fluctuations of the environment. The top Lyapunov exponent of the system  characterizes the population growth rate, and its sign determines whether the population persists or dies out.
  \medskip

The following examples illustrate the fact that the process $(\omega)_{t \geq 0}$ can be deterministic (Examples  \ref{ex:period} and \ref{ex:quasiperiod}), or stochastic (Examples \ref{ex:switch} and \ref{ex:sde}).
\bex[Periodic case]
\label{ex:period}
{\rm
Suppose $S = \RR/\ZZ$  identified with the unit circle and $$\omega_t = s + t \, ({\sf mod} \, 1)$$ for some $s \in S.$
This is the case considered in \cite{carmona2022}. Observe that here $\mu$ is the Lebesgue normalized measure on $S.$}
\eex
\bex[Quasi-periodic case]
\label{ex:quasiperiod}
{\rm A natural generalization of Example \ref{ex:period} is as follows. Suppose $S = (\RR/\ZZ)^n$ is the $n$-torus and
$$\omega_t  = (s_1 + t a_1, s_2 + t a_2, \ldots, s_n + t a_n) \, ({\sf mod} \, 1)$$ for some $s = (s_1, \ldots, s_n) \in S$ and $(a_1, \ldots, a_n)$ rationally independent numbers. That is $\sum_{i = 1}^n k_i a_i  \neq 0$ for any integers $k_1, \ldots, k_n$ such that $(k_1, \ldots, k_n) \neq (0, \ldots, 0).$ Again $(\omega_t)_{t \geq 0}$ is uniquely ergodic with  $\mu$  the Lebesgue measure on $S.$}
\eex
\bex[Switching]
\label{ex:switch}
{\rm Suppose $S = \{1, \ldots, n\}$ for some $n \in \NN^*$ and $(\omega_t)_{t \geq 0}$ is an irreducible continuous time Markov chain on $S.$
In other words, the infinitesimal generator of $(\omega_t)_{t \geq 0}$ writes $$Lf(i) = \sum_{j = 1}^n a_{ij}(f(j) - f(i))$$ for all $f : S \mapsto \RR,$ where  $(a_{ij})$ is  an irreducible rate matrix. Then $(\omega_t)_{t \geq 0}$ is uniquely ergodic and $\mu$ is the unique probability vector solution to $$\sum_{j = 1}^n (\mu_j a_{ji} - \mu_i a_{ij}) = 0$$ for all $i = 1, \ldots, n$. This situation has been considered in \cite{BS2019}.
}
\eex
\bex
\label{ex:sde}
{\rm Suppose $S$ is a compact connected Riemannian manifold and $(\omega_t)_{t \geq 0}$ a Brownian motion (or an elliptic diffusion or  more generally, the solution to a uniquely ergodic stochastic differential equation) on  $S$.  Then $(\omega_t)_{t \geq 0}$ is uniquely ergodic and $\mu$ is the normalized volume on $S$ (or a measure absolutely continuous with respect to the volume, in the diffusion case).}
\eex
We now pass to the analysis of the long term behavior of (\ref{eq:ode}).

Let $\Delta := \Delta^{d-1} = \{x \in \RR^d_+: \: \sum_{i = 1}^d x_i = 1\}$ be the unit $d-1$ simplex. Every $y \in \RR^d_+ \setminus \{0\}$ can be written as
$$y = \rho \theta,$$ with $\rho = \langle y, {\bf 1} \rangle = \sum_{i = 1}^d y_i > 0$ and $\theta = \frac{y}{\langle y, {\bf 1} \rangle} \in \Delta.$ Here and throughout, ${\bf 1}$ stands for the vector $(1, \ldots, 1)^t,$ and $\langle \cdot, \cdot \rangle$ is the usual Euclidean scalar product on $\RR^d.$

Using this decomposition, the differential equation (\ref{eq:ode}) rewrites
\beq
\label{eq:oderho}
  \frac{d \rho}{dt} = \rho \langle A(\omega_{t}) \theta, 1 \rangle
  \eeq
  and
  \beq
  \label{eq:odesimplex}
  \frac{d  \theta}{dt} = F(\omega_t, \theta),
  \eeq
  where for all $(s,\theta) \in S \times  \Delta$
  \beq
  \label{eq:defF}
 F(s,\theta) =  A(s) \theta - \langle A(s) \theta, {\bf 1} \rangle \theta.
 \eeq
   The following proposition is proved in  \cite{BS2019}, Proposition 2.13, in the case corresponding to Example \ref{ex:switch}. It mainly relies on the Random Perron-Frobenius theorem as proved by Ruelle \cite{Ruelle79} and later by   Arnold, Demetrius and Gundlach \cite{Arnold94} (see also \cite{Mier13b}, \cite{Mierc2014},  and the references therein). The proof given in \cite{BS2019}  extends to the general situation considered here. Details are given in the next section.
   \bprop
   \label{prop:main}
   Let $(\rho_t, \theta_t)$ be solution to ((\ref{eq:oderho}), (\ref{eq:odesimplex})). The
   process $(\omega_t, \theta_t)_{t \geq 0}$ is a Feller Markov process uniquely ergodic on $S  \times \Delta.$

   Let $\pi$ denote its (unique) invariant probability and let
   $$\Lambda = \int_{S \times \Delta} \langle A(s) \theta, {\bf 1} \rangle \pi(ds d \theta).$$ Then, for every initial conditions $\rho(0) > 0, \theta(0) \in \Delta$ and $\omega_0 = s,$ with probability one,
   $$\lim_{t \rar \infty} \frac{ \log (\rho_t)}{t} = \Lambda.$$
   \eprop
    For further notice, we call $\Lambda$  the {\em top Lyapunov} exponent\footnote{see the remark \ref{rem:crauel} for a justification of this terminology} of the system given by (\ref{eq:ode}).  For periodic linear differential equations
     it corresponds to what is sometimes called the principal Lyapunov exponent \cite{Mierc2014}, or  the largest Floquet multiplier \cite{carmona2022}. That is, the Floquet exponent with the largest real part. For further details we refer the reader to the Section II.2 of the excellent survey \cite{Mierc2014} by Mierczy{\'n}ski.

The following corollary easily follows from Proposition \ref{prop:main}. It provides  simple estimates of $\Lambda.$   Other estimates, mainly for periodic systems, can be found in \cite{Mierc2014} and in \cite{Mier13} for more general systems.

    \bcor The following inequalities hold true:

    \bdes
    \iti
    $$ \int_S [\min_{i = 1, \ldots, d} \sum_{j = 1}^d A_{ji}(s)] \mu(ds) \leq \Lambda \leq \int_S [\max_{i = 1, \ldots, d} \sum_{j = 1}^d A_{ji}(s)] \mu(ds);$$
    \itii
    $$\int_{S} \lambda_{min} (\frac{A(s) + A(s)^t}{2})\mu(ds) \leq  \Lambda \leq \int_{S} \lambda_{max}  (\frac{A(s) + A(s)^t}{2})\mu(ds),$$
     where $\lambda_{min}$ (respectively $\lambda_{max}$) stands for the smallest (largest) eigenvalue.
     \edes
     \ecor
     \prf  $(i).$ For all $\theta \in \Delta$ and $s \in S$ $$\min_{i = 1, \ldots, d} \sum_{j = 1}^d A_{ji}(s) \leq  \langle A(s) \theta, 1 \rangle \leq  \max_{i = 1, \ldots, d} \sum_{j = 1}^d A_{ji}(s),$$ and the result follows from the integral representation of $\Lambda.$

     $(ii).$ Let $\|y\|_{2} =  \sqrt{\langle y, y \rangle}.$ Then,  $$\frac{d}{dt} \log(\|y_t\|_2)= \frac{\langle y, A(\omega_t) y_t \rangle}{\|y_t\|_2^2} = \frac{\langle \theta_t, A(\omega_t) \theta_t \rangle}{\|\theta_t\|_2^2}.$$
    Unique ergodicity of $(\omega_t,\theta_t)_{t \geq 0},$ then implies (see e.g~ \cite{BenaimHurth}, Propositions 7.1 and  4.58) that  for every initial condition $y_0 \in \RR^d_+ \setminus \{0\},$
    $$\Lambda = \lim_{t \rar \infty} \frac{\log(\|y_t\|_2)}{t} = \int_{S \times \Delta} \frac{\langle \theta, A(s) \theta \rangle}{\|\theta\|_2^2} \pi(ds d\theta)$$ almost surely.
     Now, for all $u \in \RR^d_+$ such that $\|u\|_2 = 1,$
    $$\lambda_{min} (\frac{A(s) + A(s)^t}{2}) \leq \langle u, A(s) u \rangle \leq  \lambda_{max}  (\frac{A(s) + A(s)^t}{2}).$$ This proves the result.
    \qed
   In the particular case of a periodic system (Example \ref{ex:period}), more can be said.
   \bprop
   \label{prop:casperiodique} Suppose $S = \RR/\ZZ \backsim [0,1[$  as in Example \ref{ex:period}. There exists a continuous $1$-periodic function $t \in \RR \rar \theta^*(t) \in \Delta,$ such that:  For all $s \in S$  and  $\omega_t = s + t \, ({\sf mod} \, 1),$  $t \rar \theta^*(s+t)$ is the unique $1$-periodic solution  to (\ref{eq:odesimplex}). It is globally asymptotically stable  in the sense that    $$\lim_{t \rar \infty} \|\theta(t) - \theta^*(s+t)\| = 0$$ for every  solution $(\theta(t))_{t \geq 0}$ to  (\ref{eq:odesimplex}) with $\omega_t = s + t \, ({\sf mod} \, 1).$   In particular,
   $$\pi(ds d\theta) = ds \delta_{\theta^*(s)}(d\theta)$$ and
   $$\Lambda = \int_0^1 \langle A(s) \theta^*(s), {\bf 1} \rangle ds.$$
   \eprop
\subsection{Slow and fast regimes}
For all $T > 0,$ let $\omega^T_t = \omega_{t/T}.$ Like $(\omega_t)_{t \geq 0}$, $(\omega^T_t)_{t \geq 0}$ is a Feller Markov process on $S,$ uniquely ergodic with invariant probability $\mu.$
The parameter $1/T$ can be understood as a velocity parameter. For instance, in the context of  Example \ref{ex:period},
  $(\omega^T_t)_{t \geq 0}$ is a $T$-periodic signal. In the context of Example \ref{ex:switch},
 its  mean sojourn time in each state $i \in S$ is proportional to $T.$

  Consider the differential equation (\ref{eq:ode}) with $(\omega_t)_{t \geq 0}$ replaced by $(\omega^T_t)_{t \geq 0}.$ We let $\pi^T$ and $\Lambda^T$ denote the corresponding invariant probabilities on $S \times \Delta$ and top Lyapunov exponent as defined in Proposition \ref{prop:main}. This section considers the fast and slow regimes obtained as  $T \rar 0$ and $T \rar \infty.$
%

For a  $d \times d$ real matrix $M,$ we let $\lambda_{max}(M)$ denote the largest real part of its eigenvalues (sometimes called the {\em spectral abscissa} of $M$; see e.g~\cite{carmona2022}).

For $r > 0$ sufficiently large, $\bar{A} + r I$ has nonnegative entries and is irreducible. Hence, by Perron-Frobenius theorem (applied to  $\bar{A} + r I$), $\lambda_{max}(\bar{A})$ is an eigenvalue and there exists a unique vector, the Perron-Frobenius vector of $\bar{A},$  $\theta^* \in \Delta,$ such that $$\bar{A} \theta^* = \lambda_{max}(\bar{A}) \theta^*.$$
\bprop[Fast regime]
\label{prop:fast}
$$\lim_{T \rar 0} \pi^T = \mu \otimes \delta_{\theta^*}$$ (for the weak* topology) and
$$\lim_{T \rar 0} \Lambda^T = \lambda_{max}(\overline A).$$
\eprop
Note that Proposition \ref{prop:fast} has been proven for Example \ref{ex:switch} in (\cite{BS2019}, Corollary 2.15). The next result generalizes \cite{carmona2022}  beyond Example \ref{ex:period}. Let $\mathsf{supp}(\mu)$ be the topological support of $\mu.$ Assume that for all $s \in \mathsf{supp}(\mu),$ $A(s)$ is irreducible. Then, under this assumption, there exists for all $s \in  \mathsf{supp}(\mu)$  a unique Perron-Frobenius vector for $A(s), \theta^*(s) \in \Delta$ characterized by $$A(s) \theta^*(s) = \lambda_{max}(A(s)) \theta^*(s).$$
\bprop[Slow regime]
\label{prop:slow}
Assume that for  all $s \in \mathsf{supp}(\mu),$ $A(s)$ is irreducible.
 Then
$$\lim_{T \rar \infty} \pi^T = \mu(ds) \delta_{\theta^*(s)}$$ (for the weak* topology) and
$$\lim_{T \rar \infty} \Lambda^T = \int_S \lambda_{max}(A(s)) \mu(ds).$$
\eprop
\section{Proofs}
\label{sec:proofs}
\subsection*{Notation and Background}
 If $X$ is a  metric space (such as $S, \Delta, S \times \Delta$) we let $B(X)$ denote the space of real valued Borel bounded functions on $X$ and  $C(X) \subset B(X)$ the subspace of bounded continuous functions.
 For all $f \in B(X)$ we let $\|f\|_{\infty} = \sup_{x \in X} |f(x)|.$ If $\nu$ is a probability on $X$ and $f \in B(X)$ we write $\nu(f)$ for $\int_X f d\nu.$

 Our main assumption that $(\omega_t)_{t \geq 0}$ is a Feller Markov process on $S,$ means, as usual, that $(\omega_t)_{t \geq 0}$ is
a Markov process whose transition semigroup $(P_t)_{t \geq 0}$ is Feller. That is:
  \bdes
  \ita $P_t(C(S)) \subset C(S);$
  \itb  $\lim_{t \rar 0} P_t f(s) = f(s)$ for all $f \in C(S)$ and $s \in S.$
  \edes
It turns out (see e.g~\cite{Kallenberg}, Theorem 19.6) that $(a)$ and $(b)$ make $(P_t)_{t  \geq 0}$ strongly continuous in the sense that
$\lim_{t \rar 0} \|P_t f - f\|_{\infty} = 0$ for all $f   \in C(S).$

An {\em invariant probability} for $(\omega_t)_{t \geq 0}$ (or $(P_t)_{t  \geq 0}$) is a probability $\mu$ on $S$ such that for all $t \geq 0,$ $\mu P_t = \mu$ (i.e $\mu (P_t f) = \mu(f)$ for all $f \in B(S)$). Feller continuity and compactness of $S$ imply that such a $\mu$ always exists (see e.g~\cite{BenaimHurth}, Corollary 4.21). Our assumption  that $(\omega_t)_{t \geq 0}$ is {\em uniquely ergodic}  means that  $\mu$ is  unique.

A useful consequence of Feller continuity is that  we can assume without loss of generality that  $(\omega_t)_{t \geq 0}$ is defined on the space $\Omega$ consisting of {\em c\`adl\`ag} (right-continuous, left limit) paths $\omega: \Rp \rar S$ equipped with the Skorohod topology and associated Borel sigma field (see e.g~\cite{Kallenberg}, Theorem 19.15). As usual, for all $s \in S$ we let $\mathbb{P}_s$ denote the law of  $(\omega_t)_{t \geq 0}$ starting from $\omega_0 = s$ and  $\mathbb{P}_{\mu} = \int_S \mathbb{P}_{s} \mu(ds).$ The associated expectations are denoted $\mathbb{E}_s$ and $\mathbb{E}_{\mu}.$

For all $\omega \in \Omega$ and $t \geq 0$ we let $\mathbf{\Theta}_t(\omega)$ denote the shifted path defined as $\mathbf{\Theta}_t(\omega)(s) = \omega(t+s).$
Ergodicity of $\mu$ for the Markov process $(\omega_t)_{t \geq 0}$ makes $\mathbb{P}_{\mu}$ ergodic (but not uniquely ergodic) for the dynamical system $(\mathbf{\Theta}_t)_{t \geq 0}$ on $\Omega$ (see e.g~\cite{BenaimHurth}, Proposition 4.49).

\subsection{Proof of Propositions \ref{prop:main} and \ref{prop:casperiodique}}
For $\omega \in \Omega$ the solution to (\ref{eq:ode}) writes
$y(t) = \Phi(t,\omega) x$ where $(\Phi(t,\omega))_{t \geq 0}$ is solution to the matrix valued differential equation
$$\frac{dM}{dt} = A(\omega_t) M, M(0) = Id.$$
Let $\cal{M}_+ \subset {\cal M}$ denote the set of $d \times d$ Metzler matrices having positive diagonal entries and  $\cal{M}_{++} \subset {\cal M}_{+}$ the set of matrices having positive entries. Observe that
$$\Phi(t,\omega) \in \cal{M}_+$$ for all $t \geq 0.$ Indeed, for  $r$ large enough and all $s \in S, \, A(s) + 2r Id \geq r Id$ so that $e^{2r t} \Phi(t,\omega) \geq e^{R t} Id$ (componentwise).

For all $\theta \in \Delta$, the solution to (\ref{eq:odesimplex}) with initial condition $\theta(0) = \theta,$ writes
$$\theta(t) = \Psi(t,\omega)\theta:= \frac{\Phi(t,\omega) \theta}{\langle \Phi(t,\omega) \theta, \bf{1} \rangle}.$$

\blem
\label{lem:contraction} For $\mathbb{P}_{\mu}$ almost all $\omega \in \Omega:$
\bdes \iti There exists $N \in \NN$ such that $\Phi(t,\omega) \in \cal{M}_{++}$ for all $t \geq N;$
\itii For all $\theta, \theta' \in \Delta$ $$\lim_{t \rar \infty} \|\Psi(t,\omega) \theta - \Psi(t,\omega) \theta'\|= 0.$$
\edes
\elem
\prf
$(i).$
First observe that $\Phi(t,\omega) \in {\cal M}_{++} \Leftrightarrow  e^{rt} \Phi(t,\omega) \in {\cal M}_{++}$ for all $r > 0.$
Therefore, replacing $A(s)$ by $A(s) + r Id$ for $r > \|A\|_{\infty},$ we can assume without loss of generality that $A(s) \in {\cal M}_+$ for all $s \in S.$

Let $x(t) = \Phi(t,\omega) x$ with $x \in \RR^d_+ \setminus \{0\}.$
Suppose $x_i(0) > 0.$ Then $x_i(t) > 0$ because $\dot{x}_i(t) \geq A_{ii}(\omega_t) x_i(t) \geq 0.$
 By irreducibility of $\bar{A}$, for all $j \neq i$ there exists a sequence $i_0 = i, i_1, \ldots, i_n = j$ such that $\bar{A}_{i_k i_{k-1}} > 0$ for $k = 1,\ldots n.$
 By ergodicity, there exists a Borel set $\tilde{\Omega} \subset \Omega$ with  $\mathbb{P}_{\mu}(\tilde{\Omega}) = 1$ such that for   all $\omega \in \tilde{\Omega}$
 $$\frac{1}{t} \int_0^t A(\omega_u)du \rar \bar{A}.$$
 Therefore, for all $\omega \in \tilde{\Omega},$ there exists a sequence $t_1 > t_2 > \ldots > t_n$ with $$A_{i_k i_{k-1}}(\omega_{t_k}) > 0.$$
  By right continuity of $(\omega_t)$ we also have $A_{i_k i_{k-1}}(\omega_{t}) > 0$ for $t_k \leq t \leq t_k + \eps$ for some $\eps > 0.$
  It follows that $\dot{x}_{i_1}(t) \geq A_{i_1,i}(\omega_t) x_1(t) > 0$ for all $t_1 \leq t \leq t_1 + \eps.$
  Hence $x_{i_1}(t) > 0$ for all $t > t_1.$ Similarly $x_{i_2}(t) > 0$ for all $t > t_2$ and, by recursion, $x_j(t) > 0$ for all $t > t_n.$  In summary, we have shown that for all $i,j \in \{1, \ldots, d\}$ and $\omega \in \tilde{\Omega},$ there exists a time $t_n$ depending on $i,j, \omega$ such that for all $t \geq t_n$ $x_j(t) > 0$ whenever $x_i(0) > 0.$ This proves $(i).$

$(ii).$ Let $\RR^d_{++} = \{x  \in \RR^d \: : x_i > 0, \mbox{ for all }  i = 1\ldots d\}$ and $\dot{\Delta} = \Delta \cap \RR^d_{++}$ be the relative interior of $\Delta.$
The {\em projective} or \textit{Hilbert metric}  $d_H$ on $\RR^d_{++}$ (see Seneta \cite{Seneta}) is defined by
\[ d_H(x,y) = \log \frac{\max_{1 \leq i \leq d} x_i/y_i}{\min_{1 \leq i \leq d} x_i/y_i}. \]
Note that for all $\alpha, \beta > 0,$ $d_H(\alpha x, \beta y) = d_H(x,y)$
so that $d_H$ is not a distance on $\RR^d_{++}.$ However its restriction to $\dot{\Delta}$ is. Furthermore, for all $\theta$, $\theta' \in \dot{\Delta}$,
\begin{equation}
 \max_{1 \leq i \leq d} | \theta_i - \theta'_i| \leq e^{d_H(\theta,\theta')} - 1.
 \label{hilbnorm}
\end{equation}
 By a theorem of Birkhoff (see e.g \cite{Seneta}, Section 3.4), for all $M \in {\cal M}_{+},$
  \beq
  \label{contractdH}
   \sup_{\{x,y \in \RR^d_{++} \: d_H(x,y) > 0 \}} \frac{d_H(Mx, My)}{d_H(x,y)} = \tau[M]
   \eeq
   where $0 \leq \tau(M) \leq 1$  is the number
  defined as
$\tau(M) = \frac{1-\sqrt{r(M)}}{1+\sqrt{r(M)}}$ with $r(M) = \min_{i,j,k,l} \min \frac{M_{ik} M_{jl}}{M_{jk} M_{il}}$ if $M \in {\cal M}_{++}$
and $r(M) = 0$ if $M \in {\cal M}_{+} \setminus {\cal M}_{++}.$ In particular, for $M \in {\cal M}_+$, $\tau(M) < 1$ if and only if $M \in {\cal M}_{++}$.

  For all $0 \leq s \leq t,$ let  $$F_{s,t}(\omega) =  \max\{\log(\tau[\Phi(t-s, \mathbf{\Theta}_s(\omega))]), s-t\} \in [s-t,0].$$
 We claim that  $(F_{s,t})_{0 \leq s \leq t}$ is a sub-additive process. That is:
     \bdes
     \iti $F_{s,t} \circ \mathbf{\Theta}_v = F_{s+v,t+v},$ and
     \itii $F_{s,u} \leq F_{s,t} + F_{t,u},$
      \edes
  for all $s \leq t \leq u$ and $v \geq 0.$

    The first assertion is immediate because $\mathbf{\Theta}_s \circ \mathbf{\Theta}_v = \mathbf{\Theta}_{s+v}.$ For the second, by the cocycle property $$\Phi(u-s,\mathbf{\Theta}_s(\omega)) = \Phi(u-t,\mathbf{\Theta}_t(\omega)) \circ \Phi(t-s,\mathbf{\Theta}_s(\omega)).$$
    Thus, $$\log(\tau[\Phi(u-s,\mathbf{\Theta}_s(\omega))] \leq \log(\tau[\Phi(u-t,\mathbf{\Theta}_t(\omega))] +  \log(\tau[\Phi(t-s,\mathbf{\Theta}_s(\omega))].$$ This proves $(ii)$.

        Note also that $t,s \rar F_{s,t}(\omega)$ is continuous and that $\displaystyle \sup_{0 \leq s \leq t \leq 1} |F_{s,t}| \leq 1,$ so that the  integrability conditions  required for  the continuous time version of   Kingman's subadditive ergodic theorem (as stated in \cite{Krengel}, Theorem 5.6) are satisfied. Therefore, by this theorem,
    $$\limsup_{t \rar \infty} \frac{\log(\tau[\Phi(t,\omega)])}{t} \leq \lim_{t \rar \infty} \frac{F_{0,t}(\omega)}{t} = \gamma,$$ $\mathbb{P}_{\mu}$ almost surely, where $$\gamma =  \inf_{t > 0} \mathbb{E}_{\mu}  \frac{F_{0,t}}{t}.$$  Clearly $\gamma < 0.$ For otherwise we would have that $\tau[\Phi(n,\omega)] = 1 \Leftrightarrow \Phi(n,\omega) \in {\cal M}_+ \setminus {\cal M}_{++}$ for all $n \in \NN, \, \mathbb{P}_{\mu}$ almost surely, in contradiction with $(i).$

   Let $N$ be like in assertion $(i)$ of the Lemma. Then, by what precedes,  $\mathbb{P}_{\mu}$ almost surely,
  $$\limsup_{t \rar \infty} \frac{\log(d_H(\Psi(t+N, \omega)\theta,\Psi(t+N,\omega)\theta'))}{t}$$
   $$\leq \limsup_{t \rar \infty} \frac{\log(\tau[\Phi(t,\mathbf{\Theta}_N(\omega))])}{t} + \limsup_{t \rar \infty} \frac{\log(d_H(\Psi(N, \omega)\theta,\Psi(N,\omega)\theta'))}{t} = \gamma.$$  By inequality (\ref{hilbnorm}), this concludes the proof.
  \qed
Let $(Q_t)_{t  \geq 0}$ denote the semigroup of the process $(\omega_t,\theta_t)_{t \geq 0}.$ Then, for all $f \in B(S \times \Delta)$ $(s,\theta) \in S \times \Delta,$
$$Q_t f(s,\theta) = \mathbb{E}_s [f(\omega_t,\Psi(t,\omega)(\theta))].$$
\blem
\label{lem:feller}
The semigroup $(Q_t)_{t \geq 0}$ is Feller.
\elem
\prf We need to show that $(a)$ $Q_t(C(S \times \Delta) \subset C(S \times \Delta)$ and $(b)$ $\lim_{t \rar 0} Q_t f(s,\theta) = f(s,\theta)$ for all $f \in C(S \times \Delta).$

$(a).$ It is easy to verify that there exist constants $c_1,c_2 \geq 0$  such that for all $s,s' \in S, \theta, \theta' \in \Delta$
 \begin{eqnarray}
 \label{eq:LipF}
    \|F(s,\theta) - F(s,\theta')\| &\leq& c_1 \|\theta - \theta'\| \\
   \nonumber
    \|F(s,\theta) - F(s',\theta)\| &\leq & c_2\|A(s) - A(s')\|
 \end{eqnarray}
  where $F$ is defined by (\ref{eq:defF}).
 Fix $\eps > 0$ and let $\tilde{\omega}$ be the path defined as $\tilde{\omega}_u = \omega_{k\eps}$ for all $k\eps \leq u < (k+1) \eps.$
 Then, by  Gronwall's lemma,
 \beq
 \label{eq:gronwall}
 \|\Psi(t,\omega)(\theta) - \Psi(t,\tilde{\omega})(\theta)\| \leq c_t \int_0^t  \|A(\omega(u)) - A(\tilde{\omega}(u))\| du
 \eeq
where $c_t = e^{c_1 t} c_2.$
Thus, by Jensen inequality,
  $$\mathbb{E}_s( \|\Psi(t,\omega) \theta - \Psi(t,\tilde{\omega}))\theta)\|)^2
   \leq \mathbb{E}_s( \|\Psi(t,\omega)(\theta) - \Psi(t,\tilde{\omega})(\theta)\|^2 )$$ $$ \leq c_t^2 t \int_0^t \mathbb{E}_s ( \|A(\omega_u)) - A(\tilde{\omega}_u)\|^2) du$$
  The choice of the norm being arbitrary we  can assume that the norm on the right hand side of  the preceding inequality is the Euclidean on $\RR^{d^2}.$ Then, for all $k\eps \leq u < (k+1) \eps,$
 $$\mathbb{E}_s ( \|A(\omega_u)) - A(\tilde{\omega}_u)\|^2)  = \mathbb{E}_s\left( \mathbb{E}(\|A(\omega_u)) - A(\tilde{\omega}_u)\|^2) | {\cal F}_{k\eps})\right)$$
  $$= \mathbb{E}_s ( P_{u-k\eps}(\|A\|^2)(\omega_{k\eps}) - 2 \langle A(\omega_{k\eps}), P_{u-k\eps}(A)(\omega_{k\eps}) \rangle + \|A(\omega_{k\eps})\|^2)$$
$$ \leq \sup_{0 \leq h \leq \eps} \|P_h (\|A\|^2) - \|A\|^2)\|_{\infty} + 2 \|A\|_{\infty} \|P_h A - A \|_{\infty} := \delta(\eps).$$

   Observe that $\delta(\eps) \rar 0$ as $\eps \rar 0$ by strong continuity of $(P_t)_{t   \geq 0}.$
   Combining the two last inequalities, we get
   \beq
   \label{eq:estimateE_s}
   \mathbb{E}_s( \|\Psi(t,\omega)(\theta) - \Psi(t,\tilde{\omega})(\theta)\|^2 ) \leq c_t^2t^2 \delta(\eps)
   \eeq
   Let now $f \in C(S \times \Delta).$ Then, for every $\delta > 0$ there exists $\alpha > 0,$ such that
   $$|f(s,\theta) - f(s,\theta')| \leq \delta  + 2\|f\| {\bf 1}_{\|\theta - \theta'\| \geq \alpha}$$
   Thus
   $$|Q_t f(s,\theta) - \mathbb{E}_s(f(\omega_t, \Psi(t,\tilde{\omega})(\theta))|$$
   $$ \leq
    \mathbb{E}_s \left (|f(\omega_t, \Psi(t,\omega)(\theta)) - f(\omega_t, \Psi(t,\tilde{\omega})(\theta))|\right)
     \leq \delta + 2 \frac{\|f\|}{\alpha^2} c_t^2 t^2 \delta(\eps). $$
This shows that the left hand term goes to $0$ uniformly in $(s,\theta) \in S \times \Delta$ as $\eps \rar 0.$

In order to conclude it  suffices  to show that $(s,\theta) \rar \mathbb{E}_s \left (f(\omega_t, \Psi(t,\tilde{\omega})(\theta)) \right ) $ is continuous.
 For all $s \in S,$ let $(\Psi^s_t)_{t \geq 0}$ denote the semi-flow on $\Delta$ induced by the autonomous differential equation
 $\frac{d\theta}{dt} =  F(s,\theta)$
 Then for $k\eps \leq t < (k+1) \eps$
 $$f(\omega_t, \Psi(t,\tilde{\omega})(\theta)) = f(\omega_t, \Psi^{\omega_{k \eps}}_{t-k\eps} \circ ... \Psi^{\omega_{\eps}}_{\eps} \circ \Psi^{\omega_0}_{\eps}(\theta))$$
Now, for  every $h \in C(S^{k+2} \times \Delta),$   Feller continuity of $(P_t)_{t  \geq 0}$, makes the map $(s,\theta) \rar \mathbb{E}_s(h(\omega_t, \omega_{k\eps}, \ldots, \omega_0, \theta)$  continuous. This is immediate to verify when $h$ is a product function (i.e~ $h(s_{k+1}, \ldots, s_0,\theta) = h_{k+1}(s_{k+1}) \cdot  h_0(s_0) g(\theta)$) and the general case follows by the density in $C(S^{k+2} \times \Delta)$ of the vector space span by  product functions. This concludes the proof of $(a).$

$(b).$ Let $f \in C(S  \times \Delta)$ and $\delta > 0.$
Because $\|\Psi(t,\omega)(\theta) - \theta\| \leq t \|F\|_{\infty}$,
$\|Q_t f(s,\theta) - \mathbb{E}_s (f(\omega_t,\theta))\| \leq \delta$ for all $t$ sufficiently small. By Feller continuity of $(P_t)$ $\lim_{t \rar 0} \mathbb{E}_s (f(\omega_t,\theta)) = \lim_{t \rar 0} P_t (f(\cdot, \theta))(s) = f(s,\theta).$
\qed
We can now conclude the proof of Proposition \ref{prop:main}.
It follows from Lemma \ref{lem:contraction} $(ii)$ that for all $f : S \times \Delta \mapsto \RR$ continuous (hence uniformly continuous) and all $\theta, \theta' \in \Delta,$ $$\lim_{t \rar \infty} |f(\omega_t,\Psi(t,\omega) \theta) - f(\omega_t,\Psi(t,\omega) \theta') |  = 0$$ $\mathbb{P}_s$ almost surely, for $\mu$ almost all $s \in S.$
Hence, for all $\theta, \theta' \in \Delta,$
\beq
\label{eq:contractQ}
\lim_{t \rar \infty} |Q_t f(s,\theta) - Q_t f(s,\theta')| = 0,
\eeq for $\mu$ almost all $s \in S.$
Let now $\pi$ be an invariant probability of $(Q_t)_{t \geq 0}.$ Such a $\pi$ always exist because $(Q_t)_{t \geq 0}$ if Feller on $S \times \Delta$ compact. To prove that $\pi$ is unique, assume that $\pi'$ is another invariant probability.
Then, writing $\pi(f)$ for $\int_{S \times \Delta} f(s,\theta) \pi(ds d\theta),$
$$\pi(f) - \pi'(f) = \pi (Q_t f) - \pi' (Q_t f)$$
 $$= \int_S \left [ \int_{\Delta \times \Delta} ( Q_t f(s,\theta) - Q_t f(s,\theta') ) \pi(d \theta|s) \pi(d \theta'|s) \right ]\mu(ds)$$
 where for each $s \in S,$ $\pi(.|s)$ (respectively $\pi'(.|s)$)  is a conditional distribution of $\pi$ (respectively $\pi'$)  (see \cite{Dud02}, Section 10.2). It then follows from (\ref{eq:contractQ}) and dominated convergence that $\pi(f) =   \pi'(f).$ Thus $\pi = \pi'.$ This proves unique ergodicity.

 Now, unique ergodicity and Feller continuity of $(\omega_s,\theta_s)_{s \geq 0}$ imply that  for every continuous function $g : S \times \Delta \rar \RR$ $$\lim_{t \rar \infty} \frac{1}{t} \int_0^t g(\omega_s,\theta_s) ds = \int g d\pi$$ $\mathbb{P}_{s,\theta}$ almost surely for all $s,\theta \in S \times \Delta.$ (see e.g~ \cite{BenaimHurth}, Proposition 7.1 for discrete time chains combined with Proposition 4.58 to handle continuous time).
 This concludes the proof of Proposition \ref{prop:main} with $g(s,\theta) = \langle A(s) \theta, {\bf 1} \rangle)$.
 \brem
 {\rm
 \label{rem:crauel}
  By the multiplicative ergodic theorem, there exist numbers $\Lambda_1  < \ldots < \Lambda_r, r \leq d,$ called {\em Lyapunov exponents}, such that for $\mathbb{P}_{\mu}$ almost all $\omega$ and all $x \in \RR^d \setminus \{0\},$  $$\lim_{t \rar \infty} \frac{\log{\|\Phi(t,\omega)x\|}}{t} := \Lambda(x,\omega) \in \{\Lambda_1, \ldots, \Lambda_r\}.$$ The set of $x \in \RR^d$ for which $\Lambda(x,\omega) < \Lambda_r$ is a  vector space (depending on $\omega$) having nonzero codimension. On the other hand, by what precedes, $\Lambda(x,\omega) = \Lambda$ for all $x \in \RR^d_+ \setminus \{0\}.$ It follows that $\Lambda = \Lambda_r.$
 }
 \erem
 \subsection*{Proof of Proposition \ref{prop:casperiodique}}
 For $s \in S = \RR/\ZZ,$  let $\omega[s] \in \Omega$ be the path defined as  $$\omega_t[s] = s + t \, ({\sf mod} \, 1).$$
 By Brouwer fixed point theorem, the map $\Psi(1,\omega[0]):\Delta \mapsto \Delta$ has a fixed point $\theta^*.$ Set $\theta^*(t) =\Psi(t,\omega[0])(\theta^*).$ Then $$\theta^*(t+1) = \Psi(t, \omega[1]) \circ \Psi(1,\omega[0])(\theta^*) = \theta^*(t)$$ proving that $t \rar \theta^*(t)$ is $1$-periodic.

 For all $s \in S$ and $\theta \in \Delta$
 $$\lim_{t \rar \infty} \|\Psi(t,\omega[s])(\theta)-\theta^*(t+s)\| = \lim_{t \rar \infty} \|\Psi(t,\omega[s])(\theta)- \Psi(t,\omega[s])(\theta^*(s))\| = 0,$$
by Lemma \ref{lem:contraction} applied with $\omega = \omega[s]$ and $\theta' = \theta^*[s].$ Observe here, that  the conclusions of Lemma \ref{lem:contraction} hold with $\omega = \omega[s]$ for all $s \in S$ simply because  $\omega[s]$ is $1$-periodic.
\brem
\mylabel{rem:floquetcomparison}
{\rm  The proof given here can be re-interpreted in the classical framework of Floquet's theory used in \cite{carmona2022}.  By Floquet's theorem, every solution to $\frac{dy}{dt} = A(\omega_t[0]) y$  writes $y(t) = P(t) e^{t B} y(0),$ where $P$ is a $1$-periodic matrix such that $P(0) = Id.$ The matrix $e^B$ has nonnegative entries, hence, by Perron-Frobenius theorem, a  eigenvector $y^* \in \RR^d_{+} \setminus \{0\}.$ The point $\theta^*$ in the proof above is the projection of $y^*$ on the simplex, $\theta^* = y^*/|y^*|.$}
\erem
 \section{Proof of Propositions \ref{prop:fast} and \ref{prop:slow}}
   For all $T > 0,$ let $(Q^T_t)_{t  \geq 0}$ denote the semigroup of $(\omega_t^T,\theta_t)_{t \geq 0}$ with $\omega_t^T = \omega_{t/T}$ and $(\theta_t)_{t \geq 0}$ is  solution to (\ref{eq:odesimplex}) when $\omega_t$ is replaced by $\omega^T_t.$
 Using the notation of the preceding section one sees that
 $$Q^T_t (f)(s,\theta) = \mathbb{E}_s \left[f(\omega_{t/T}, \Psi(t, \omega^T)(\theta))\right]$$ for all $f \in B(S \times \Delta).$
 \subsection*{Proof of Proposition \ref{prop:fast}}
 For all $\theta \in \Delta,$ let $\bar{F}(\theta) = \int_S F(s,\theta) \mu(ds),$ where $F$ is defined by (\ref{eq:defF}). Let $(\bar{\Psi}_t)_{t \geq 0}$ denote the semi-flow on $\Delta$ induced by the
 differential equation $\displaystyle \dot{\theta} = \bar{F}(\theta).$ The following lemma follows from the averaging principle  as given in Freidlin and Wentzell \cite{FW12}(Theorem 2.1, Chapter 7).
 \blem
 \label{lem:averaging}
 For all $\delta > 0$ and  $t \geq 0,$
 $$\lim_{T \rar 0} \mathbb{P}_{\mu} \left (\sup_{\theta \in \Delta, 0 \leq u \leq t} \|\Psi(u, \omega^T)(\theta) - \bar{\Psi}_u(\theta)\| \geq \delta \right) = 0.$$ In particular, for all $f \in C(\Delta)$ and $t \geq 0,$
 $$\lim_{T \rar 0} \mathbb{E}_{\mu} \left [\|f \circ \Psi(t,\omega^T) - f \circ \bar{\Psi}_t\|_{\infty} \right] = 0.$$
 \elem
 \prf
 We claim that
 \beq
 \label{eq:condaveraging}
 \lim_{R \rar \infty}  \sup_{t \geq 0} \mathbb{P}_{\mu} \left( \left |\frac{1}{R} \int_t^{t+R} F(\omega_s,\theta) ds - \bar{F}(\theta)\right | \geq \delta \right ) = 0.
  \eeq
 Indeed, by stationarity (invariance of $\mathbb{P}_{\mu}$ for $(\mathbf{\Theta}_t)_{t \geq 0}$),
 $$\mathbb{P}_{\mu} \left( \left |\frac{1}{R} \int_t^{t+R} F(\omega_s,\theta) ds - \bar{F}(\theta)\right | \geq \delta \right ) = \mathbb{P}_{\mu} \left( \left |\frac{1}{R} \int_0^{R} F(\omega_s,\theta) ds - \bar{F}(\theta)\right | \geq \delta \right )$$  for all $t \geq 0;$ and the right hand term goes to $0,$ as $R \rar \infty,$ by ergodicity of $\mu.$

  By the averaging theorem  (Theorem 2.1, Chapter 7 in \cite{FW12}), condition (\ref{eq:condaveraging}) implies  that for all $\delta > 0, t \geq 0$ and $\theta \in \Delta,$
 \beq
 \label{eq:averaging}
 \lim_{T \rar 0} \mathbb{P}_{\mu} \left (\sup_{0 \leq u \leq t} \|\Psi(u, \omega^T)(\theta) - \bar{\Psi}_u(\theta)\| \geq \delta \right) = 0.
 \eeq
By Lipschitz continuity  (see (\ref{eq:LipF})) and Gronwall's  lemma,
  $$\sup_{0 \leq u \leq t} \|\Psi(u, \omega^T)(\theta) - \Psi(u,\omega^T)(\theta')\| +  \|\bar{\Psi}_t(\theta) - \bar{\Psi}_t(\theta')\|\leq 2 e^{c_1t}\|\theta -\theta'\|$$ for all $\theta, \theta' \in \Delta.$ Fix  $\eps < \frac{\delta}{4} e^{-c_1t}$ and let  $\{B(\theta_i, \eps), i = 1, \ldots, N\}$ be a finite covering of $\Delta$ by balls of radius $\eps.$ Then
 $$\sup_{0 \leq u \leq t, \theta \in \Delta} \|\Psi(u, \omega^T)(\theta) - \bar{\Psi}_u(\theta)\| \leq \max_{i=1, \ldots, N} \sup_{0 \leq u \leq t} \|\Psi(u, \omega^T)(\theta_i) - \bar{\Psi}_u(\theta_i)\| + \delta/2.$$
 Hence
 $$\mathbb{P}_{\mu} \left (\sup_{0 \leq u \leq t, \theta \in \Delta} \|\Psi(u, \omega^T)(\theta) - \bar{\Psi}_u(\theta)\| \geq \delta \right) \leq \sum_{i = 1}^N \mathbb{P}_{\mu} \left (\sup_{0 \leq u \leq t} \|\Psi(u, \omega^T)(\theta_i) - \bar{\Psi}_u(\theta_i)\| \geq \delta/2\right).$$ The right hand term goes $0$ as $T \rar 0$ by (\ref{eq:averaging}).
  \qed
 We now prove the proposition. Let $\pi^T$ be the invariant measure of $(Q^T_t)_{t \geq 0}$ and let $\pi^{0}$ be a limit point of $(\pi^T)_{T > 0}$ for the weak* topology, as $T \rar 0.$
  That is: $\pi^{T_n} f \rar \pi^{0} f$ for some sequence $T_n \rar 0$ and all $f \in C(S \times \Delta).$

 Let  $p : S \times \Delta \rar \Delta$ be the projection defined as $p(s,\theta) = \theta$ and let  $\pi_2^T = \pi^T \circ p^{-1}$ be the second marginal of $\pi^T.$ Similarly, set $\pi^0_2 = \pi^{0} \circ p^{-1}.$

 For all $f \in C(\Delta)$ and $t \geq 0,$
 $$\pi_2^T f = \pi^T (f \circ p) = \pi^T Q_t^T (f \circ p) = \int_{S \times \Delta} \mathbb{E}_s[f(\Psi(t,\omega^T)(\theta))]\pi^T(ds d\theta).$$
 Thus,
 $$|\pi_2^T f -\pi_2^T (f \circ \bar{\Psi}_t)| = \left|\int_{S \times \Delta} \mathbb{E}_s[f(\Psi(t,\omega^T)(\theta)) -f(\bar{\Psi}_t(\theta))] \pi^T(ds d\theta)\right|$$
 $$\leq \int_S \mathbb{E}_s [ \|f\circ \Psi(t,\omega^T) -f \circ \bar{\Psi}_t\|_{\infty}] \mu(ds) = \mathbb{E}_{\mu} [ \|f\circ \Psi(t,\omega^T) -f \circ \bar{\Psi}_t\|_{\infty}].$$
 Here we have used the fact that the first marginal of $\pi^T$ is $\mu.$
Using Lemma \ref{lem:averaging}, it comes that
 $$\pi_2^0 f = \pi_2^0(f \circ \bar{\Psi}_t)$$ for all $t \geq 0.$ This proves that $\pi_2^0$ is invariant for  $\{\bar{\Psi}_t\}_{t \geq 0},$ but since $\{\bar{\Psi}_t\}_{t \geq 0}$ has $\theta^*$ as globally asymptotically stable equilibrium,  necessarily $\pi_2^{0} =\delta_{\theta^*}.$ On the other hand, the first marginal of $\pi^0$ is $\mu.$ Thus $\pi^0 = \mu \otimes \delta_{\theta^*}.$  This concludes the proof.
 \subsection*{Proof of Proposition \ref{prop:slow}}
 Recall (see the proof of Lemma \ref{lem:feller}) that for all $s \in S,$ we let $(\Psi^s_t)_{t \geq 0}$ denote the semi-flow on $\Delta$ induced by the
 differential equation $\displaystyle \dot{\theta} = F(s,\theta).$

 Let $(Q^{\infty}_t)_{t \geq 0})$ denote the Markov semigroup  on $S \times \Delta$
defined as    $$Q^{\infty}_t f(s,\theta) = f(s,\Psi^s_t(\theta))$$ for all $f \in B(S \times \Delta).$
 \blem
 \label{lem:approxlargeT}
  For all $f \in C(S \times \Delta)$ and  $t \geq 0$
$$\lim_{T \rar \infty} \|Q^T_t f - Q^{\infty}_t f\|_{\infty} = 0.$$
 \elem
 \prf Let  $f \in C(S \times \Delta).$ By uniform continuity of $f,$ for  every $t > 0$ and $\delta > 0$ there exists $\alpha > 0$ such that
$$|f(\omega^T_t,\Psi(t,\omega^T) \theta) - f(s,\Psi^s_t(\theta))| \leq \delta + 2\|f\|_{\infty}
 {\bf 1}_{\large \{ d(\omega^T_t, s) + \|\Psi(t,\omega^T) (\theta) - \Psi^s_t(\theta)\| \geq \alpha \}}.$$
Thus
$$\|Q^T_t f(s,\theta) - Q^{\infty}_t f(s,\theta) \| \leq \mathbb{E}_s \left( |f(\omega^T_t,\Psi(t,\omega^T) (\theta) - f(s,\Psi^s_t(\theta))| \right)$$
$$\leq \delta + 2\|f\|_{\infty} \frac{ \mathbb{E}_s (\|\Psi(t,\omega^T) (\theta) - \Psi^s_t(\theta)\|) + P_{t/T} (d(.,s))(s) }{\alpha}.$$
By Feller continuity, $P_{t/T}(d(.,s))(s) \rar 0$ uniformly in $s \in S$ as $T \rar \infty.$ This follows for example from Lemma 19.3 ($F_3$) in \cite{Kallenberg}.
Now the estimate (\ref{eq:estimateE_s}) applied with $\omega^T$ in place of $\omega,$ $(P_{t/T})_{t \geq 0}$ in place of
 $(P_t)_{t \geq 0}$  and $\eps > t$ gives
 \beq
 \label{eq:gronwallestimate}
 \sup_{s \in S}\mathbb{E}_s (\|\Psi(t,\omega^T) (\theta) - \Psi^s_t(\theta)\|) \leq c^2 t^2 \delta(\eps/T).
 \eeq
 with $\delta(\eps/T) \rar 0$ as $T \rar \infty.$
 This concludes the proof.
 \qed
 We can now prove Proposition \ref{prop:slow}.
 Let $\pi^T$ be the invariant measure of $(\omega^T_t, \theta_t)$ and let $\pi^{\infty}$ be a limit point of $(\pi^T)_{T > 0}$ for the weak* topology. That is $\pi^{T_n} f \rar \pi^{\infty} f$ for some sequence $T_n \rar \infty$ and all $f \in C(S \times \Delta).$
 Then,
 $$|\pi^T(f) - \pi^T (Q_t^{\infty} f)|= |\pi^T (Q_t^T(f) - Q_t^{\infty}(f))| \leq \|Q_t^T(f) - Q_t^{\infty}(f)\|_{\infty}.$$
 Thus, by Lemma \ref{lem:approxlargeT}, (ii),
 $\pi^{\infty}(f) = \pi^{\infty}(Q_t^{\infty}(f)).$
 Now for all $s \in \mathsf{supp}(\mu)$ $$\lim_{t \rar \infty} Q_t^{\infty}(f)(s,\theta) = \lim_{t \rar \infty} f(s, \Psi^s_t(\theta)) = f(s, \theta^*(s)).$$ Thus, since $\pi^{\infty}(\mathsf{supp}(\mu) \times \Delta) = 1,$ it comes that
 $$\pi^{\infty}(f) = \int_S \int_{\Delta} f(s,\theta^*(s)) \pi^{\infty}(ds d\theta) = \int_S  f(s,\theta^*(s)) \mu(ds).$$
This proves the first part of Proposition \ref{prop:slow}. The second part follows directly from the first one.
\section{Concluding remarks}
The results and proofs given here all rely on the  assumption that $(\omega_t)_{t  \geq 0}$ is a Markov process.  In particular, they do not apply to the case where $t \rar \omega_t$ is a deterministic periodic signal  with discontinuities. This situation is investigated in the preprint \cite{BLSS23b}. The recent preprint \cite{MS23} provides a  first order expansion of $\Lambda^T$ when $T$ goes to $0$.
\subsection*{Acknowledgment}
This research is supported by  the Swiss National Foundation  grants 200020 196999 and 200020 219913. We thank Philippe Carmona and an anonymous referee for their careful reading and useful comments.
\bibliographystyle{abbrv}
\bibliography{LyapounovMetzler}

\begin{thebibliography}{10}

\bibitem{Arnold94}
L.~Arnold, V.~M. Gundlach, and L.~Demetrius.
\newblock Evolutionary formalism for products of positive random matrices.
\newblock {\em Ann. Appl. Probab.}, 4(3):859--901, 1994.

\bibitem{BenaimHurth}
M.~Bena\"{\i}m and T.~Hurth.
\newblock {\em Markov Chains on Metric Spaces, A Short Course}, volume~99 of
  {\em Universitext}.
\newblock Springer, Cham, 2022.

\bibitem{BLSS23a}
M.~Bena\"{\i}m, C.~Lobry, T.~Sari, and E.~Strickler.
\newblock Untangling the role of temporal and spatial variations in persistence
  of populations.
\newblock {\em Theoretical Population Biology}, 154:1--26, 2023.

\bibitem{BLSS23b}
M.~Bena\"{\i}m, C.~Lobry, T.~Sari, and E.~Strickler.
\newblock {When can a population spreading across sink habitats persist ?}
\newblock working paper, https://hal.inrae.fr/hal-04099082, May 2023.

\bibitem{BS2019}
M.~Bena\"{\i}m and E.~Strickler.
\newblock Random switching between vector fields having a common zero.
\newblock {\em Ann. Appl. Probab.}, 29(1):326--375, 2019.

\bibitem{carmona2022}
P.~Carmona.
\newblock Asymptotic of the largest {Floquet} multiplier for cooperative
  matrices.
\newblock {\em Annales de la Facult\'e des sciences de
  Toulouse:\,Math\'ematiques}, Ser. 6, 31(4):1213--1221, 2022.

\bibitem{Dud02}
R.~M. Dudley.
\newblock {\em Real analysis and probability}, volume~74 of {\em Cambridge
  Studies in Advanced Mathematics}.
\newblock Cambridge University Press, Cambridge, 2002.
\newblock Revised reprint of the 1989 original.

\bibitem{FW12}
M.~Freidlin and A.~Wentzell.
\newblock {\em Random perturbations of dynamical systems}, volume 260 of {\em
  Grundlehren der Mathematischen Wissenschaften [Fundamental Principles of
  Mathematical Sciences]}.
\newblock Springer, Heidelberg, third edition, 2012.
\newblock Translated from the 1979 Russian original by Joseph Sz\"ucs.

\bibitem{HirschSmith06}
M.~Hirsch and H.~Smith.
\newblock Chapter 4 monotone dynamical systems.
\newblock volume~2 of {\em Handbook of Differential Equations: Ordinary
  Differential Equations}, pages 239--357. North-Holland, 2006.

\bibitem{Kallenberg}
O.~Kallenberg.
\newblock {\em Foundations of modern probability}, volume~99 of {\em
  Probability Theory and Stochastic Modelling}.
\newblock Springer, Cham, 2021.
\newblock Third edition of the 1997 original.

\bibitem{Krengel}
U.~Krengel.
\newblock {\em Ergodic Theorems}, volume~99 of {\em De Gruyter Series in
  Mathematics}.
\newblock De Gruyter, 1985.

\bibitem{MalikSmith08}
T.~Malik and H.~Smith.
\newblock Does dormancy increase fitness of bacterial populations in
  time-varying environments?
\newblock {\em Bulletin of mathematical biology}, 70:1140--62, 06 2008.

\bibitem{Mierc2014}
J.~Mierczy{\'n}ski.
\newblock Estimates for principal lyapunov exponents: A survey.
\newblock {\em Nonautonomous Dynamical Systems}, 1(1), 2014.

\bibitem{Mier13}
J.~Mierczy{\'n}ski.
\newblock Lower estimates of top {L}yapunov exponent for cooperative random
  systems of linear {ODE}s.
\newblock {\em Proc. Amer. Math. Soc.}, 143(3):1127--1135, 2015.

\bibitem{Mier13b}
J.~Mierczy\'{n}ski and W.~Shen.
\newblock Principal {L}yapunov exponents and principal {F}loquet spaces of
  positive random dynamical systems. {II}. {F}inite-dimensional systems.
\newblock {\em J. Math. Anal. Appl.}, 404(2):438--458, 2013.

\bibitem{MS23}
P.~Monmarché and E.~Strickler.
\newblock Asymptotic expansion of the invariant measurefor markov-modulated
  odes at high frequency, 2023.
\newblock arXiv 2309.16464.

\bibitem{Ruelle79}
D.~Ruelle.
\newblock Analyticity properties of the characteristic exponents of random
  matrix.
\newblock {\em Adv. Math.}, 32(3):68--80, 1979.

\bibitem{Seneta}
E.~Seneta.
\newblock {\em Non-negative matrices and {M}arkov chains}.
\newblock Springer Series in Statistics. Springer, New York, 2006.
\newblock Revised reprint of the second (1981) edition.

\end{thebibliography}
\end{document}